\documentclass[thmsa,12pt]{article}
\usepackage{amssymb}

\usepackage{sw20lart}

\input{tcilatex}
\begin{document}

\begin{center}
\bigskip

$G$\textbf{\ Method\ in\ Action:}

\bigskip

\textbf{Pivot}$^{\text{+}}$\textbf{\ Algorithm for\ Self-avoiding\ Walk}

\bigskip

\textbf{Udrea\ P\u {a}un}

\bigskip
\end{center}

\noindent The pivot algorithm --- we also call it the pivot chain --- is an
algorithm for approximately uniform sampling from $\Omega _{N},$ the set of $%
N$-step self-avoiding walks on $\Bbb{Z}^{d}$ ($N,$ $d\geq 1$). Based on this
algorithm and the $G$ method, we construct another algorithm/chain, called
the pivot$^{\text{+}}$ algorithm/chain, for approximately uniform sampling
from $\Omega _{N},$ here, $N\geq 2$. The pivot$^{\text{+}}$ algorithm
samples the pivot from the set $\left\{ 1,2,...,N-1\right\} $ according to
the uniform probability distribution on this set while the pivot algorithm
samples the pivot from the set $\left\{ 0,1,2,...,N-1\right\} $ according to
the uniform probability distribution on this set, so, on the pivot, the pivot%
$^{\text{+}}$ algorithm is better than the pivot algorithm. Further, we
obtain another important thing, namely, the pivot$^{\text{+}}$
algorithm/chain enters, at time $1$, a set $2d$ times smaller than $\Omega
_{N},$ and stays forever in this set, so, at times $1,2,...$ we work with a
chain having a state space $2d$ times smaller than $\Omega _{N}.$ As to the
speed of convergence, we conjecture that the pivot$^{\text{+}}$
algorithm/chain is faster than the pivot algorithm/chain.

\medskip



\begin{center}
\medskip

\textbf{1. }$G$\textbf{\ METHOD}

\medskip
\end{center}

In this section, we present the $G$ method --- notation, notions, results.

\smallskip

Set 
\[
\text{Par}\left( E\right) =\left\{ \Delta \left| \text{ }\Delta \text{ is a
partition of }E\right. \right\} , 
\]

\noindent where $E$ is a nonempty set. We shall agree that the partitions do
not contain the empty set. $\left( E\right) $ is the improper (degenerate)
partition of $E$ --- we need this partition.

\smallskip

\textit{Definition}\textbf{\ }1.1. Let $\Delta _{1},\Delta _{2}\in $Par$%
\left( E\right) .$ We say that $\Delta _{1}$\textit{\ is finer than }$\Delta
_{2}$ if $\forall V\in \Delta _{1},$ $\exists W\in \Delta _{2}$ such that $%
V\subseteq W.$

Write $\Delta _{1}\preceq \Delta _{2}$ when $\Delta _{1}$ is finer than $%
\Delta _{2}.$

\smallskip

In this article, a vector is a row vector and a stochastic matrix is a row
stochastic matrix.

The entry $\left( i,j\right) $ of a matrix $Z$ will be denoted $Z_{ij}$ or,
if confusion can arise, $Z_{i\rightarrow j}.$

\smallskip

Set 
\[
\left\langle m\right\rangle =\left\{ 1,2,...,m\right\} \text{ (}m\geq 1\text{%
)}, 
\]

\[
N_{m,n}=\left\{ P\left| \text{ }P\text{ is a nonnegative }m\times n\text{
matrix}\right. \right\} \text{,} 
\]
\[
S_{m,n}=\left\{ P\left| \text{ }P\text{ is a stochastic }m\times n\text{
matrix}\right. \right\} \text{,} 
\]

\[
N_{n}=N_{n,n}\text{,} 
\]

\[
S_{n}=S_{n,n}. 
\]

Let $P=\left( P_{ij}\right) \in N_{m,n}$. Let $\emptyset \neq U\subseteq
\left\langle m\right\rangle $ and $\emptyset \neq V\subseteq \left\langle
n\right\rangle $. Set the matrices 
\[
P_{U}=\left( P_{ij}\right) _{i\in U,j\in \left\langle n\right\rangle },\text{
}P^{V}=\left( P_{ij}\right) _{i\in \left\langle m\right\rangle ,j\in V},%
\text{ and }P_{U}^{V}=\left( P_{ij}\right) _{i\in U,j\in V}. 
\]

\smallskip

Let $W$ be a nonempty finite set. Suppose that $W=\left\{
s_{1},s_{2},...,s_{t}\right\} $. Set 
\[
\left( \left\{ i\right\} \right) _{i\in W}\in \text{Par}\left( W\right) ,%
\text{ }\left( \left\{ i\right\} \right) _{i\in W}=\left( \left\{
s_{1}\right\} ,\left\{ s_{2}\right\} ,...,\left\{ s_{t}\right\} \right) . 
\]

\noindent \textit{E.g.}, 
\[
\left( \left\{ i\right\} \right) _{i\in \left\langle n\right\rangle }=\left(
\left\{ 1\right\} ,\left\{ 2\right\} ,...,\left\{ n\right\} \right) . 
\]

\smallskip

\textit{Definition} 1.2. Let $P\in N_{m,n}$. We say that $P$ is a \textit{%
generalized stochastic matrix} if $\exists a\geq 0$, $\exists Q\in S_{m,n}$
such that $P=aQ$.

\smallskip

\textit{Definition} 1.3 ([13]). Let $P\in N_{m,n}$. Let $\Delta \in $Par$%
\left( \left\langle m\right\rangle \right) $ and $\Sigma \in $Par$\left(
\left\langle n\right\rangle \right) $. We say that $P$ is a $\left[ \Delta
\right] $\textit{-stable matrix on }$\Sigma $ if $P_{K}^{L}$ is a
generalized stochastic matrix, $\forall K\in \Delta ,$ $\forall L\in \Sigma $%
. In particular, a $\left[ \Delta \right] $-stable matrix on $\left( \left\{
i\right\} \right) _{i\in \left\langle n\right\rangle }$ is called $\left[
\Delta \right] $\textit{-stable }for short.

\smallskip

\textit{Definition} 1.4 ([13]). Let $P\in N_{m,n}$. Let $\Delta \in $Par$%
\left( \left\langle m\right\rangle \right) $ and $\Sigma \in $Par$\left(
\left\langle n\right\rangle \right) $. We say that $P$ is a $\Delta $\textit{%
-stable matrix on }$\Sigma $ if $\Delta $ is the least fine partition for
which $P$ is a $\left[ \Delta \right] $-stable matrix on $\Sigma $. In
particular, a $\Delta $-stable matrix on $\left( \left\{ i\right\} \right)
_{i\in \left\langle n\right\rangle }$ is called $\Delta $\textit{-stable}
while a $\left( \left\langle m\right\rangle \right) $-stable matrix on $%
\Sigma $ is called \textit{stable on }$\Sigma $ for short. A stable matrix
on $\left( \left\{ i\right\} \right) _{i\in \left\langle n\right\rangle }$
is called \textit{stable} for short.

\smallskip

Let $\Delta \in $Par$\left( \left\langle m\right\rangle \right) $ and $%
\Sigma \in $Par$\left( \left\langle n\right\rangle \right) $. Set (see [13]
for $G_{\Delta ,\Sigma }$ and [14] for $\overline{G}_{\Delta ,\Sigma }$) 
\[
G_{\Delta ,\Sigma }=G_{\Delta ,\Sigma }\left( m,n\right) =\left\{ P\left| 
\text{ }P\in S_{m,n}\text{ and }P\text{ is a }\left[ \Delta \right] \text{%
-stable matrix on }\Sigma \right. \right\} 
\]

\noindent and 
\[
\overline{G}_{\Delta ,\Sigma }=\overline{G}_{\Delta ,\Sigma }\left(
m,n\right) =\left\{ P\left| \text{ }P\in N_{m,n}\text{ and }P\text{ is a }%
\left[ \Delta \right] \text{-stable matrix on }\Sigma \right. \right\} . 
\]

\smallskip

When we study or even when we construct products of nonnegative matrices (in
particular, products of stochastic matrices) using $G_{\Delta ,\Sigma }$ or $%
\overline{G}_{\Delta ,\Sigma },$ we shall refer this as the $G$\textit{\
method}. $G$ comes from the verb \textit{to group} and its derivatives.

\smallskip

The $G$\ method together or not together with the Gibbs sampler in a
generalized sense led to good or very good results --- for the first case
(``together''), see, \textit{e.g.}, [16]-[17], and, for the second one, see, 
\textit{e.g.}, [13, Examples 2.11 and 2.19] and [14, Theorems 1.6 and 1.8].
Moreover, the $G$\ method suggested another method, a $G$-type method,
namely, the $G^{+}$\ method --- for this powerful method, see [18].

\smallskip

Let $P\in \overline{G}_{\Delta ,\Sigma }.$ Then (see Definition 1.3) $%
\forall K\in \Delta ,$ $\forall L\in \Sigma ,$ $\exists a_{K,L}\geq 0,$ $%
\exists Q_{K,L}\in S_{\left| K\right| ,\left| L\right| }$ such that $%
P_{K}^{L}=a_{K,L}Q_{K,L}.$ Set (see, \textit{e.g.}, [15]) 
\[
P^{-+}=\left( P_{KL}^{-+}\right) _{K\in \Delta ,L\in \Sigma },\text{ }%
P_{KL}^{-+}=a_{K,L},\text{ }\forall K\in \Delta ,\text{ }\forall L\in \Sigma 
\]

\noindent ($P_{KL}^{-+},$ $K\in \Delta ,$ $L\in \Sigma ,$ are the entries of
matrix $P^{-+}$). If confusion can arise, we write $P^{-+\left( \Delta
,\Sigma \right) }$ instead of $P^{-+}.$ Obviously, $\forall K\in \Delta ,$ $%
\forall L\in \Sigma $ we have 
\[
a_{K,L}=\sum\limits_{j\in L}P_{ij},\text{ }\forall i\in K. 
\]

\smallskip

Let $P\in N_{m,n}.$ Let $\Delta \in $Par$\left( \left\langle m\right\rangle
\right) .$ Set 
\[
\overline{\alpha }\left( P\right) =\frac{1}{2}\max\limits_{1\leq i,j\leq
m}\sum\limits_{k=1}^{n}\left| P_{ik}-P_{jk}\right| 
\]

\noindent (see, \textit{e.g.}, [4, p. 56], [5, p. 144], and [19, p. 137])
and, more generally, 
\[
\overline{\gamma }_{\Delta }\left( P\right) =\max_{K\in \Delta }\overline{%
\alpha }\left( P_{K}\right) 
\]
\noindent (see [11]; see, \textit{e.g.}, also [12]-[13]). $\overline{\alpha }
$ and $\overline{\gamma }_{\Delta }$ are called \textit{ergodicity
coefficients} when $P$ is a stochastic matrix.

\smallskip

\textbf{THEOREM\ 1.5} ([12]-[13]). \textit{If }$P_{1}\in G_{\Delta
_{1},\Delta _{2}}\subseteq S_{m_{1},m_{2}},$ $P_{2}\in G_{\Delta _{2},\Delta
_{3}}\subseteq S_{m_{2},m_{3}},...,$ $P_{n-1}\in G_{\Delta _{n-1},\Delta
_{n}}\subseteq S_{m_{n-1},m_{n}},$ $P_{n}\in S_{m_{n},m_{n+1}},$ \textit{then%
} 
\[
\overline{\gamma }_{\Delta _{1}}\left( P_{1}P_{2}...P_{n}\right) \leq 
\overline{\gamma }_{\Delta _{1}}\left( P_{1}\right) \overline{\gamma }%
_{\Delta _{2}}\left( P_{2}\right) ...\overline{\gamma }_{\Delta _{n}}\left(
P_{n}\right) . 
\]

\smallskip

\textit{Proof.} See [12, Theorem 1.18] and [13, Remark 2.21(b)]. $\square $

\smallskip

For applications of $\overline{\gamma }_{\Delta },$ see, \textit{e.g.}, [11,
Theorems 2.18 and 3.6] and [13, Example 2.22]. Below we give another
application --- an important one --- of $\overline{\gamma }_{\Delta }.$

\smallskip

\textbf{THEOREM\ 1.6} ([14]). \textit{Let }$P_{1}\in \overline{G}_{\Delta
_{1},\Delta _{2}}\subseteq N_{m_{1},m_{2}},$ $P_{2}\in \overline{G}_{\Delta
_{2},\Delta _{3}}\subseteq N_{m_{2},m_{3}},...,P_{n}\in \overline{G}_{\Delta
_{n},\Delta _{n+1}}\subseteq N_{m_{n},m_{n+1}}.$\textit{\ Suppose that }$%
\Delta _{1}=\left( \left\langle m_{1}\right\rangle \right) $ \textit{and} $%
\Delta _{n+1}=\left( \left\{ i\right\} \right) _{i\in \left\langle
m_{n+1}\right\rangle }.$ \textit{Then} 
\[
P_{1}P_{2}...P_{n} 
\]
\noindent \textit{is a stable matrix} (\textit{i.e., a matrix with identical
rows, see Definition} 1.4). \textit{Moreover, we have} 
\[
\left( P_{1}P_{2}...P_{n}\right) _{\left\{ i\right\}
}=P_{1}^{-+}P_{2}^{-+}...P_{n}^{-+},\text{ }\forall i\in \left\langle
m_{1}\right\rangle 
\]
\noindent ($\left( P_{1}P_{2}...P_{n}\right) _{\left\{ i\right\} }$ \textit{%
is the row }$i$\textit{\ of} $P_{1}P_{2}...P_{n}$)\textit{.}

\smallskip

\textit{Proof.} See [14]. When $P_{1},$ $P_{2},...,P_{n}$ are stochastic
matrices, we can give another proof for the first part. Indeed, in this
case, by Theorem 1.5 we have 
\[
\overline{\gamma }_{\Delta _{1}}\left( P_{1}P_{2}...P_{n}\right) \leq 
\overline{\gamma }_{\Delta _{1}}\left( P_{1}\right) \overline{\gamma }%
_{\Delta _{2}}\left( P_{2}\right) ...\overline{\gamma }_{\Delta _{n}}\left(
P_{n}\right) . 
\]
\noindent Since $\Delta _{1}=\left( \left\langle m_{1}\right\rangle \right)
, $ we have 
\[
\overline{\gamma }_{\Delta _{1}}=\overline{\alpha }. 
\]
\noindent Since $P_{n}\in G_{\Delta _{n},\Delta _{n+1}}$ and $\Delta
_{n+1}=\left( \left\{ i\right\} \right) _{i\in \left\langle
m_{n+1}\right\rangle },$ we have 
\[
\overline{\gamma }_{\Delta _{n}}\left( P_{n}\right) =0. 
\]
\noindent It follows that 
\[
\overline{\alpha }\left( P_{1}P_{2}...P_{n}\right) =0, 
\]
\noindent and, further, that 
\[
P_{1}P_{2}...P_{n} 
\]
\noindent is a stable matrix. $\square $

\smallskip

Below we consider a similarity relation, and give a result about it --- the
relation and result are, in this article, useful for the construction of
pivot$^{\text{+}}$ algorithm/chain, see Section 2, Comment 1 to Theorem 2.8
and its proof.

\smallskip

\textit{Definition} 1.7 ([15]). Let $P,Q\in \overline{G}_{\Delta ,\Sigma
}\subseteq N_{m,n}.$ We say that $P$\textit{\ is similar to }$Q$ (\textit{%
with respect to }$\Delta $\textit{\ and }$\Sigma $) if 
\[
P^{-+}=Q^{-+}. 
\]

Set $P\backsim Q$ when $P$ is similar to $Q.$ Obviously, $\backsim $ is an
equivalence relation on $\overline{G}_{\Delta ,\Sigma }.$

\smallskip

Set $e=e\left( n\right) =\left( 1,1,...,1\right) \in \Bbb{R}^{n},$ $\forall
n\geq 1$; $e^{\prime }$ is the transpose of $e$.

\smallskip

\textit{Example} 1.8. Let 
\[
P=\left( 
\begin{array}{cccc}
\frac{1}{4}\smallskip  & \frac{1}{4} & \frac{1}{4} & \frac{1}{4} \\ 
\frac{1}{4} & \frac{1}{4}\smallskip  & \frac{1}{4} & \frac{1}{4} \\ 
\frac{1}{4} & \frac{1}{4} & \frac{1}{4}\smallskip  & \frac{1}{4} \\ 
\frac{1}{4} & \frac{1}{4} & \frac{1}{4} & \frac{1}{4}\smallskip 
\end{array}
\right) ,\text{ }Q=\left( 
\begin{array}{cccc}
\frac{2}{4}\smallskip  & 0 & \frac{2}{4} & 0 \\ 
\frac{2}{4} & 0\smallskip  & \frac{2}{4} & 0 \\ 
\frac{2}{4} & 0 & \frac{2}{4}\smallskip  & 0 \\ 
\frac{2}{4} & 0 & \frac{2}{4} & 0\smallskip 
\end{array}
\right) ,\text{ and }R=\left( 
\begin{array}{cccc}
\frac{1}{4}\smallskip  & \frac{1}{4} & \frac{2}{4} & 0 \\ 
\frac{1}{4} & \frac{1}{4}\smallskip  & \frac{2}{4} & 0 \\ 
\frac{2}{4} & 0 & 0\smallskip  & \frac{2}{4} \\ 
0 & \frac{2}{4} & \frac{2}{4} & 0\smallskip 
\end{array}
\right) .
\]
\noindent We have $P,Q,R\in G_{\Delta ,\Sigma }\subseteq S_{4},$ where $%
\Delta =\left( \left\langle 4\right\rangle \right) $ and $\Sigma =\left(
\left\{ 1,2\right\} ,\left\{ 3,4\right\} \right) ,$ and 
\[
P^{-+}=Q^{-+}=R^{-+}=\left( 
\begin{array}{cc}
\frac{2}{4} & \frac{2}{4}
\end{array}
\right) .
\]
\noindent So, $P\backsim Q\backsim R.$ Note that: 1) $P$ is a stable
(stochastic) matrix, \textit{i.e.}, a matrix with identical rows, and $%
P=e^{\prime }\delta ,$ where $\delta =\left( \frac{1}{4},\frac{1}{4},\frac{1%
}{4},\frac{1}{4}\right) $ --- obviously, $e=\left( 1,1,1,1\right) $; 2) any
matrix which is similar to $P$ with respect to $\Delta $ and $\Sigma $ has
at least $8$ and at most $16$ positive entries; 3) $Q$ has the smallest
possible number of positive entries, \textit{i.e.}, $8,$ while $P$ and $R$
have not --- $P$ has the greatest possible number of positive entries, 
\textit{i.e.}, $16$.

\smallskip

\textbf{THEOREM\ 1.9} ([15]). \textit{Let }$P_{1},U_{1}\in \overline{G}%
_{\Delta _{1},\Delta _{2}}\subseteq N_{m_{1},m_{2}},$\textit{\ }$%
P_{2},U_{2}\in \overline{G}_{\Delta _{2},\Delta _{3}}\subseteq
N_{m_{2},m_{3}},$\textit{\ }$...,$\textit{\ }$P_{n},U_{n}\in \overline{G}%
_{\Delta _{n},\Delta _{n+1}}\subseteq N_{m_{n},m_{n+1}}.$\textit{\ Suppose
that} 
\[
P_{1}\backsim U_{1},\mathit{\ }P_{2}\backsim U_{2},\mathit{\ }%
...,P_{n}\backsim U_{n}. 
\]
\noindent \textit{Then} 
\[
P_{1}P_{2}...P_{n}\backsim U_{1}U_{2}....U_{n}. 
\]

\noindent \textit{If, moreover, }$\Delta _{1}=\left( \left\langle
m_{1}\right\rangle \right) $\textit{\ and }$\Delta _{n+1}=\left( \left\{
i\right\} \right) _{i\in \left\langle m_{n+1}\right\rangle },$\textit{\ then}
\[
P_{1}P_{2}...P_{n}=U_{1}U_{2}....U_{n} 
\]

\noindent (\textit{therefore, when }$\Delta _{1}=\left( \left\langle
m_{1}\right\rangle \right) $\textit{\ and }$\Delta _{n+1}=\left( \left\{
i\right\} \right) _{i\in \left\langle m_{n+1}\right\rangle },$\textit{\ a
product of }$n$\textit{\ representatives, the first of an equivalence class
included in }$\overline{G}_{\Delta _{1},\Delta _{2}},$\textit{\ the second
of an equivalence class included in }$\overline{G}_{\Delta _{2},\Delta
_{3}}, $\textit{\ }$...,$\textit{\ the }$n$th\textit{\ of an equivalence
class included in }$\overline{G}_{\Delta _{n},\Delta _{n+1}},$\textit{\ does
not depend on the choice of representatives --- we can work }(\textit{this
is important})\textit{\ with }$U_{1}U_{2}....U_{n}$ \textit{instead of} $%
P_{1}P_{2}...P_{n}$ \textit{and}, \textit{conversely, with} $%
P_{1}P_{2}...P_{n}$ \textit{instead of} $U_{1}U_{2}....U_{n}$; \textit{due
to Theorem} 1.6, $P_{1}P_{2}...P_{n}$ \textit{and} $U_{1}U_{2}....U_{n}$ 
\textit{are stable matrices}).

\smallskip

\textit{Proof}. See [15]. $\square $

\smallskip

\textit{Example} 1.10. Consider the matrices $Q$ and $R$ from Example 1.8.
Let 
\[
P_{1}=Q,\text{ }U_{1}=R,\text{ and }P_{2}=U_{2}=\left( 
\begin{array}{cccc}
\frac{1}{3}\smallskip  & \frac{2}{3} & 0 & 0 \\ 
\frac{1}{3} & \frac{2}{3}\smallskip  & 0 & 0 \\ 
0 & 0 & \frac{2}{5}\smallskip  & \frac{3}{5} \\ 
0 & 0 & \frac{2}{5} & \frac{3}{5}
\end{array}
\right) .
\]
\noindent We have $P_{1},$ $U_{1}\in G_{\Delta _{1},\Delta _{2}}\subseteq
S_{4}$ and $P_{2},$ $U_{2}\in G_{\Delta _{2},\Delta _{3}}\subseteq S_{4},$
where
\[
\Delta _{1}=\left( \left\langle 4\right\rangle \right) ,\text{ }\Delta
_{2}=\left( \left\{ 1,2\right\} ,\left\{ 3,4\right\} \right) ,\text{ }\Delta
_{3}=\left( \left\{ i\right\} \right) _{i\in \left\langle 4\right\rangle }
\]
\noindent ($\Delta _{3}=\left( \left\{ 1\right\} ,\left\{ 2\right\} ,\left\{
3\right\} ,\left\{ 4\right\} \right) $). Since $P_{1},$ $U_{1}\in G_{\Delta
_{1},\Delta _{2}}\subseteq S_{4},$ $P_{2},$ $U_{2}\in G_{\Delta _{2},\Delta
_{3}}\subseteq S_{4},$ $P_{1}\backsim U_{1},\mathit{\ }P_{2}\backsim U_{2},$ 
$\Delta _{1}=\left( \left\langle 4\right\rangle \right) ,$ and $\Delta
_{3}=\left( \left\{ i\right\} \right) _{i\in \left\langle 4\right\rangle },$
by Theorem 1.9 we have $P_{1}P_{2}=U_{1}U_{2}$ --- this thing can also be
obtained by direct computation; we have 
\[
P_{1}P_{2}=U_{1}U_{2}=\left( 
\begin{array}{cccc}
\frac{2}{12}\smallskip  & \frac{4}{12} & \frac{4}{20} & \frac{6}{20} \\ 
\frac{2}{12} & \frac{4}{12}\smallskip  & \frac{4}{20} & \frac{6}{20} \\ 
\frac{2}{12} & \frac{4}{12} & \frac{4}{20}\smallskip  & \frac{6}{20} \\ 
\frac{2}{12} & \frac{4}{12} & \frac{4}{20} & \frac{6}{20}
\end{array}
\right) =e^{\prime }\tau ,
\]

\noindent where $\tau =\left( \frac{2}{12},\frac{4}{12},\frac{4}{20},\frac{6%
}{20}\right) $ ($P_{1}P_{2}$ and $U_{1}U_{2}$ are stable matrices).

\bigskip

\begin{center}
\textbf{2. PIVOT\ AND PIVOT}$^{\text{+}}$\textbf{\ ALGORITHMS}
\end{center}

\bigskip

In this section, we present two algorithms, the pivot algorithm (known) and
pivot$^{\text{+}}$ algorithm (new), for approximately uniform sampling ---
unfortunately, not for uniform sampling --- from the set of $N$-step
self-avoiding walks. For these algorithms, certain fundamental results are
given. Finally, we conjecture that the pivot$^{\text{+}}$ algorithm is
faster than the pivot algorithm.

\smallskip

Set

\[
\left\langle \left\langle m\right\rangle \right\rangle =\left\{
0,1,...,m\right\} \text{ (}m\geq 0\text{)}. 
\]

\noindent Set 
\[
\left| \left| x\right| \right| _{2}=\sqrt{\sum\limits_{i=1}^{n}x_{i}^{2}}, 
\]

\noindent where $x=\left( x_{1},x_{2},...,x_{n}\right) \in \Bbb{R}^{n}$ and $%
n\geq 1.$ $\left| \left| x\right| \right| _{2}$ is called the 2-\textit{norm}
or \textit{Euclidean norm of }$x$; $\left| \left| \cdot \right| \right| _{2}$
is called the 2-\textit{norm} or \textit{Euclidean norm.}

\smallskip

\textit{Definition} 2.1 ([2]; see, \textit{e.g.}, also [7] and [9]). Let $d,$
$N\geq 1.$ Let $\omega =\left( \omega _{0},\omega _{1},...,\omega
_{N}\right) ,$ where $\omega _{i}\in \Bbb{Z}^{d},\forall i\in \left\langle
\left\langle N\right\rangle \right\rangle .$ We say that $\omega $ is a ($N$-%
\textit{step}) \textit{self-avoiding walk }(\textit{on} $\Bbb{Z}^{d}$) if $%
\omega _{0}=0,$ $\left| \left| \omega _{i+1}-\omega _{i}\right| \right|
_{2}=1,$ $\forall i\in \left\langle \left\langle N-1\right\rangle
\right\rangle ,$ and $\omega _{i}\neq \omega _{j},$ $\forall i,j\in
\left\langle \left\langle N\right\rangle \right\rangle ,$ $i\neq j.$ $N$ is
called the \textit{length of} (\textit{self-avoiding walk}) $\omega .$

\smallskip

The self-avoiding walks on $\Bbb{Z}^{d}$ were invented by the chemist Paul
Flory [2] to study the polymers.

\smallskip

The self-avoiding walks on $\Bbb{Z}^{d}$ can be represented by the
(nondirected infinite) graph, called lattice or grid, and having the vertex
set $\Bbb{Z}^{d}$ and the edge set 
\[
\left\{ \left\{ x,y\right\} \left| \text{ }x,y\in \Bbb{Z}^{d}\text{ and }%
\left| \left| x-y\right| \right| _{2}=1\right. \right\} 
\]
\noindent --- the edge set is a set of 2-element subsets of (the vertex set) 
$\Bbb{Z}^{d}$. A self-avoiding walk on $\Bbb{Z}^{d}$ can be seen as being a
walk on the above graph, starting from $0$ (from the origin), and which does
not visit any vertex more than once. (See, \textit{e.g.}, [1]--[3], [7],
[9], and [20].)

\smallskip

Set 
\[
\Omega _{N}=\left\{ \omega \left| \omega \text{ is an }N\text{-step
self-avoiding walk on }\Bbb{Z}^{d}\right. \right\} . 
\]

\smallskip

Denote by $\mathcal{O}_{d}$ the set of orthogonal (linear) transformations
of $\Bbb{R}^{d}$ which leave the above lattice/grid invariant (each such a
transformation leaves the origin fixed (because it is a linear
transformation) and maps the lattice into itself). (See, \textit{e.g.}, [9,
p. 323] and, for $d=2,$ [7, p. 51].)

\smallskip

The next algorithm is due to M. Lal [6], and is called the \textit{pivot
algorithm} --- see, \textit{e.g.}, also [9, pp. 314 and 323] and, for $d=2,$
[7, p. 51].

\smallskip

\textit{Step} 0. Start with a ($N$-step) self-avoiding walk, say, $\omega .$

\smallskip

\textit{Step} 1. Choose a number, called \textit{pivot} in our article, from 
$\left\langle \left\langle N-1\right\rangle \right\rangle $ according to the
uniform distribution on this set. Suppose that this pivot is, say, $k.$

\smallskip

\textit{Step} 2. Choose a transformation from $\mathcal{O}_{d}$ according to
the uniform distribution on this set ($\mathcal{O}_{d}$ is a finite set).
Suppose that this transformation is, say, $T.$

\smallskip

\textit{Step} 3. Consider the walk obtained by fixing the first $k$ steps of
(the self-avoiding walk) $\omega ,$ but performing the transformation $T$ on
the remaining part of $\omega ,$ using $\omega _{k}$ as the origin for the
transformation --- we obtain a walk which may or may not be self-avoiding ($%
\omega \in \Omega _{N}$ fixed, $i\longmapsto \omega _{i},$ $i\in
\left\langle \left\langle N-1\right\rangle \right\rangle ,$ is a bijective
function, so, we can call either $k$ or $\omega _{k}$ pivot --- in [20, p.
66] $\omega _{k}$ is called the \textit{pivot site}...).

\smallskip

\textit{Step} 4. If the walk obtained at Step 3 is self-avoiding, choose it
as the new walk, and go to Step 1. If the walk obtained at Step 3 is not
self-avoiding, remove/reject it and keep $\omega ,$ and go to Step 1.

\smallskip

The above algorithm is, in fact, a (homogeneous finite) Markov chain, so, we
will use ``the pivot algorithm'' or ``the pivot chain'' to denote it. Denote
by $P$ the transition matrix of this chain and by $q_{n}$ the probability
distribution at time $n$ of this chain, $\forall n\geq 0.$ $q_{0}$ is also
called the \textit{initial probability distribution} (of this chain). Step 0
(of the pivot algorithm) says that 
\[
\left( q_{0}\right) _{\alpha }=\left\{ 
\begin{array}{cc}
1 & \text{if }\alpha =\omega , \\ 
0 & \text{if }\alpha \neq \omega ,
\end{array}
\right. 
\]
\noindent $\forall \alpha \in \Omega _{N}$ ($q_{0}=\left( \left(
q_{0}\right) _{\alpha }\right) _{\alpha \in \Omega _{N}}$). Recall that
(see, \textit{e.g.}, [5, p. 137] and [19, p. 114]) 
\[
q_{n+1}=q_{n}P=q_{0}P^{n},\text{ }\forall n\geq 0. 
\]

\smallskip

Let $M\in \left\langle N\right\rangle .$ Let $\gamma =\left( \gamma
_{0},\gamma _{1},...,\gamma _{M}\right) \in \Omega _{M}.$ Set 
\[
K_{\gamma }=K_{\left( \gamma _{0},\gamma _{1},...,\gamma _{M}\right)
}=\left\{ \omega \left| \text{ }\omega \in \Omega _{N}\text{ and }\omega
_{i}=\gamma _{i},\forall i\in \left\langle \left\langle M\right\rangle
\right\rangle \right. \right\} . 
\]

Set 
\[
E_{d}=\left\{ e_{1},e_{2},...,e_{d},-e_{1},-e_{2},...,-e_{d}\right\} , 
\]

\noindent where $e_{1},e_{2},...,e_{d}$ are the (row) vectors of canonical
basis of $\Bbb{R}^{d}.$

\smallskip

\textit{Remark} 2.2. $\gamma _{1}\in E_{d}$ because $\forall N\geq 1,$ $%
\forall \omega \in \Omega _{N}$ we have 
\[
\omega _{1}\in E_{d}. 
\]

\smallskip

\textbf{THEOREM 2.3.} (i) (this result is known, see, \textit{e.g.}, [9, p.
324]) $P$ \textit{is irreducible}.

\smallskip

(ii) $P_{K_{\left( \gamma _{0},\gamma _{1}\right) }}^{K_{\left( \gamma
_{0},\gamma _{1}\right) }}$ \textit{is irreducible} (\textit{nonnegative,
substochastic}), $\forall \gamma _{1}\in E_{d}$ ($\gamma _{0}=0;$ \textit{%
for }$P_{K_{\left( \gamma _{0},\gamma _{1}\right) }}^{K_{\left( \gamma
_{0},\gamma _{1}\right) }},$\textit{\ see Section }1 (\textit{see the
definition of matrix} $P_{U}^{V}...$)).

\medskip

\textit{Proof}. (i) See, \textit{e.g.}, [9, pp. 350$-$353] --- the proof
from there has two cases, Case I and Case II, and is based on the fact that
any self-avoiding walk from $\Omega _{N}$ can be unwrapped to form a
(completely) straight self-avoiding walk. ($P$ is irreducible if and only if 
$\forall \alpha ,\beta \in \Omega _{N},$ either $P_{\alpha \beta }>0$ or $%
\exists t\geq 1,$ $\exists \nu ^{\left( 1\right) },\nu ^{\left( 2\right)
},...,\nu ^{\left( t\right) }\in \Omega _{N}$ such that $P_{\alpha \nu
^{\left( 1\right) }},P_{\nu ^{\left( 1\right) }\nu ^{\left( 2\right)
}},...,P_{\nu ^{\left( t\right) }\beta }>0,$ see, \textit{e.g.}, [4, p. 52],
[9, p. 286], and [19, pp. 11 and 18].)

\smallskip

(ii) The proof is based on the fact that any self-avoiding walk from $%
K_{\left( \gamma _{0},\gamma _{1}\right) }$ can be unwrapped to form a
straight self-avoiding walk from $K_{\left( \gamma _{0},\gamma _{1}\right) }$
($K_{\left( \gamma _{0},\gamma _{1}\right) }$ contains just one straight
self-avoiding walk), all the self-avoiding walks of this unwrapping process
being from $K_{\left( \gamma _{0},\gamma _{1}\right) }.$ The proof of this
fact is simple, see, \textit{e.g.}, [9, pp. 350$-$353], see, there, Case I
--- $t$ from there is greater or equal to $1$ --- and Case II --- $q$ from
there is greater or equal to $1$; $t\geq 1$ and $q\geq 1$ imply the above
fact. $\square $

\smallskip

\textit{Problem} 2.4. Fix $\tau =\left( \tau _{0},\tau _{1},...,\tau
_{N}\right) \in \Omega _{N}.$ Fix $M_{0}\in \left\langle N\right\rangle .$
Find the smallest number $M\in \left\langle N\right\rangle ,$ $M\geq M_{0},$
such that
\[
P_{K_{\left( \tau _{0},\tau _{1},...,\tau _{M}\right) }}^{K_{\left( \tau
_{0},\tau _{1},...,\tau _{M}\right) }}
\]
\noindent is irreducible ($M\leq N$ when $M_{0}\in \left\langle
N-1\right\rangle ;$ $M=M_{0}$ when $M_{0}=1$ (Theorem 2.3(ii)) or when $%
M_{0}=N$).

\smallskip

Denote by $\pi $ the uniform probability distribution on $\Omega _{N}$.
Recall that $e=e\left( n\right) =\left( 1,1,...,1\right) \in \Bbb{R}^{n},$ $%
\forall n\geq 1$; $e^{\prime }$ is the transpose of $e$.

\smallskip

Besides the fact that the transition matrix $P$ of pivot chain is
irreducible (Theorem 2.3(i)), this chain has other important proprieties ---
see the next result.

\smallskip

\textbf{THEOREM 2.5.} (Known things.) (i) $P$ \textit{is aperiodic
irreducible and symmetric}.

(ii) $\pi P=\pi .$

(iii) $P^{n}\longrightarrow e^{\prime }\pi $ \textit{as} $n\longrightarrow
\infty $ ($e^{\prime }\pi $ \textit{is a stable matrix, i.e., a matrix with
identical rows, see Definition} 1.4).

(iv) $q_{n}\longrightarrow \pi $ \textit{as} $n\longrightarrow \infty $ (%
\textit{therefore}, $\pi $ \textit{is the limit probability distribution of
pivot chain}).

\smallskip

\textit{Proof}. (i) $P$ is aperiodic irreducible because it is irreducible
(Theorem 2.3(i)) and $P_{\alpha \alpha }>0,$ $\forall \alpha \in \Omega
_{N}. $ $P$ is symmetric (\textit{i.e.}, $P_{\alpha \beta }=P_{\beta \alpha
},$ $\forall \alpha ,\beta \in \Omega _{N}$) because any orthogonal
transformation, being bijective, is invertible.

(ii) $\pi P=\pi $ follows from $\pi _{\alpha }P_{\alpha \beta }=\pi _{\beta
}P_{\beta \alpha },$ $\forall \alpha ,\beta \in \Omega _{N}.$

(iii) See, \textit{e.g.}, [4, p. 123, Theorem 4.2] (in [4], a vector is, see
p. 46, a column vector), [8, pp. 54$-$55], and [9, pp. 286$-$287].

(iv) $q_{n}=q_{0}P^{n}\longrightarrow q_{0}e^{\prime }\pi =\pi $ as $%
n\longrightarrow \infty .$ $\square $

\smallskip

We cannot run the pivot algorithm/chain forever --- we stop it at a fixed
time, say, $s;$ $s$ must be as large as possible, $q_{s}$ (the probability
distribution at time $s$ (of pivot chain)) being (see Theorem 2.5(iv)) close
or not too close (nobody knows --- Markov chain Monte Carlo!) to $\pi $ (the
uniform probability distribution on $\Omega _{N}$). Due to this fact, we
sample (exact sampling) from $\Omega _{N}$ according to $q_{s}$, not
according to $\pi .$

\smallskip

\textbf{THEOREM 2.6. }$\left( K_{\left( \gamma _{0},\gamma _{1},...,\gamma
_{M}\right) }\right) _{\left( \gamma _{0},\gamma _{1},...,\gamma _{M}\right)
\in \Omega _{M}}$ \textit{is a partition of }$\Omega _{N},$\textit{\ }$%
\forall M\in \left\langle N\right\rangle .$\textit{\ If, moreover, }$M=1,$%
\textit{\ we have} ($\left| \cdot \right| $ \textit{is the cardinality}) 
\[
\left| K_{\left( \gamma _{0},\gamma _{1}\right) }\right| =\left| K_{\left(
\tau _{0},\tau _{1}\right) }\right| ,\text{ }\forall \gamma _{1},\tau
_{1}\in E_{d}
\]
(\textit{equivalently, }$\forall \left( \gamma _{0},\gamma _{1}\right)
,\left( \tau _{0},\tau _{1}\right) \in \Omega _{1}$\textit{\ --- see Remark}
2.2 \textit{again}); \textit{further, we have --- recall that the
probability distribution }$\pi $\textit{\ }(\textit{on }$\Omega _{N}$)%
\textit{\ is} \textit{uniform}; $P=P\left( \cdot \right) $ \textit{is the
probability} --- 
\[
P\left( K_{\left( \gamma _{0},\gamma _{1}\right) }\right) =\frac{1}{2d},%
\text{ }\forall \gamma _{1}\in E_{d}.
\]

\smallskip

\textit{Proof}. The first part is obvious. Further, we suppose that\textit{\ 
}$M=1$. Fix $\gamma _{1}$ and $\tau _{1}$ ($\gamma _{0}=\tau _{0}=0$).
Suppose that $\gamma _{1}\neq \tau _{1}$. Consider the rotation, say, $R$
with (the property that) $R\left( \gamma _{1}\right) =\tau _{1}.$ (More
generally, we can consider an orthogonal (linear) transformation, say, $T$
with $T\left( \gamma _{1}\right) =\tau _{1}.$) Consider the function 
\[
F:K_{\left( \gamma _{0},\gamma _{1}\right) }\longrightarrow K_{\left( \tau
_{0},\tau _{1}\right) },\text{ }F\left( \omega _{0},\omega _{1},...,\omega
_{N}\right) =\left( R\left( \omega _{0}\right) ,R\left( \omega _{1}\right)
,...,R\left( \omega _{N}\right) \right) . 
\]
\noindent $F$ is well-defined: $R\left( \omega _{0}\right) =R\left( 0\right)
=0=\tau _{0},$ $R\left( \omega _{1}\right) =R\left( \gamma _{1}\right) =\tau
_{1},$ and $\left( R\left( \omega _{0}\right) ,R\left( \omega _{1}\right)
,...,R\left( \omega _{N}\right) \right) $ is self-avoiding because $R$ is
bijective, so, 
\[
\left( R\left( \omega _{0}\right) ,R\left( \omega _{1}\right) ,...,R\left(
\omega _{N}\right) \right) \in K_{\left( \tau _{0},\tau _{1}\right) }. 
\]
\noindent $F$ is bijective (injective and surjective) because $R$ is
bijective (invertible...). So,

\[
\left| K_{\left( \gamma _{0},\gamma _{1}\right) }\right| =\left| K_{\left(
\tau _{0},\tau _{1}\right) }\right| ; 
\]
\noindent further, we obtain 
\[
P\left( K_{\left( \gamma _{0},\gamma _{1}\right) }\right) =\frac{1}{2d}.%
\text{ }\square 
\]

\smallskip

Set (recall that $\left| \cdot \right| $ is the cardinality) 
\[
c_{N}=\left| \Omega _{N}\right| ,\text{ }\forall N\geq 1\text{ (}N\in \Bbb{N}%
\text{)}
\]
\noindent and 
\[
a_{N}=\left| K_{\left( 0,e_{1}\right) }\right| 
\]
\noindent (by Theorem 2.6 we have 
\[
a_{N}=\left| K_{\left( \gamma _{0},\gamma _{1}\right) }\right| ,\text{ }%
\forall \gamma _{1}\in E_{d}\text{).}
\]

\smallskip

Let $A$ and $B$ be finite sets. Let $n\in \Bbb{N}$, $n\geq 1.$ We say that $%
A $ \textit{is smaller than }$B$ if $\left| A\right| \leq \left| B\right| .$
We say that $A$ \textit{is larger than }$B$ if $\left| A\right| \geq \left|
B\right| $ (equivalently, if $B$ is smaller than $A$). We say that $A$ 
\textit{is }$n$ \textit{times smaller than }$B$ if $\left| A\right| =\frac{%
\left| B\right| }{n}.$ We say that $A$ \textit{is }$n$ \textit{times larger
than }$B$ if $\left| A\right| =n\left| B\right| $ (equivalently, if $B$ is%
\textit{\ }$n$ times smaller than $A$).

\smallskip

\textbf{THEOREM 2.7.} \textit{We have} 
\[
c_{N}=2da_{N},\text{ }\forall N\geq 1 
\]
\noindent (\textit{therefore, }$\Omega _{N}$\textit{\ is }$2d$\textit{\
times larger than any set }$K_{\left( \gamma _{0},\gamma _{1}\right) },$%
\textit{\ where} $\gamma _{1}\in E_{d}$).

\smallskip

\textit{Proof}. By Theorem 2.6, 
\[
c_{N}=\left| \Omega _{N}\right| =\left| \bigcup\limits_{\gamma _{1}\in
E_{d}}K_{\left( \gamma _{0},\gamma _{1}\right) }\right| =\sum\limits_{\gamma
_{1}\in E_{d}}\left| K_{\left( \gamma _{0},\gamma _{1}\right) }\right|
=\left| E_{d}\right| a_{N}=2da_{N},\text{ }\forall N\geq 1.\text{ }\square 
\]

One way to sample exactly from $\Omega _{N}$ according to $\pi $ (the
uniform probability distribution on $\Omega _{N}$) is by means of our
reference method from [15]. Below we present this method for the above case.

\smallskip

\textit{Step} I. Choose a set from the partition $\left( K_{\left( \gamma
_{0},\gamma _{1}\right) }\right) _{\gamma _{1}\in E_{d}}$ according to the
uniform probability distribution on this partition. Suppose that the result
of sampling is, say, $K_{\left( \tau _{0},\tau _{1}\right) }.$

\smallskip

\textit{Step} II. Choose a self-avoiding walk from the set $K_{\left( \tau
_{0},\tau _{1}\right) }$ according to the uniform probability distribution
on this set. Suppose that the result of sampling is, say, $\omega $.

\smallskip

$\omega $ is the result of sampling from $\Omega _{N}$ according to $\pi .$

\smallskip

We can associate this reference method with 2 steps, Steps I and II, with a
Markov chain, a reference chain, a chain which can do what this reference
method does, see, in this section, Comment 1 to Theorem 2.8 and its proof,
see also [15] for the general case for the reference method and its
associate (Markov) chain, the reference chain (this chain was constructed
using the $G$ method) --- we do not present this method and this chain here
because their presentation is too long.

\smallskip

Denote by $B_{N}$ the set of $N$-step straight self-avoiding walks.
Obviously, 
\[
B_{N}\subseteq \Omega _{N}\text{ and }\left| B_{N}\right| =\left|
E_{d}\right| =\left| \left( K_{\left( \gamma _{0},\gamma _{1}\right)
}\right) _{\gamma _{1}\in E_{d}}\right| =2d.
\]
\noindent We will use this set to replace Step I with a step equivalent to
it ---the step obtained is due to the first matrix of the reference chain
(from [15]). Step II, unfortunately, cannot be performed when $K_{\left(
\tau _{0},\tau _{1}\right) }$ is too large ($\left| K_{\left( \tau _{0},\tau
_{1}\right) }\right| =a_{N}=\frac{\left| \Omega _{N}\right| }{2d}$). So, we
will resort to the pivot algorithm as to this step.

\smallskip

We arrived to the following algorithm for approximately uniform sampling
from $\Omega _{N}.$ This algorithm works when $N\geq 2$ --- suppose that $%
N\geq 2.$

\smallskip

\textit{Step} 0. Start with a (straight or not) self-avoiding walk, say, $%
\omega $.

\smallskip

\textit{Step} 1. Choose a straight self-avoiding walk from $B_{N}$ according
to the uniform probability distribution on this set. Suppose that this
straight self-avoiding walk is, say, $\tau .$ Remove $\omega $ and keep $%
\tau .$ (This step is equivalent to Step I --- $K_{\left( \gamma _{0},\gamma
_{1}\right) }$ contains just one straight self-avoiding walk, $\forall
\gamma _{1}\in E_{d}$...; further, see Theorem 2.8...)

\smallskip

\textit{Step} 2. Choose a number, called \textit{pivot}, from $\left\langle
N-1\right\rangle $ (warning! not from $\left\langle \left\langle
N-1\right\rangle \right\rangle $) according to the uniform probability
distribution on this set. Suppose that this pivot is, say, $k$.

\smallskip

\textit{Step} 3. Choose a transformation from $\mathcal{O}_{d}$ according to
the uniform probability distribution on this set. Suppose that this
transformation is, say, $T.$

\smallskip

\textit{Step} 4. Consider the walk obtained by fixing the first $k$ steps of 
$\tau ,$ but performing the transformation $T$ on the remaining part of $%
\tau ,$ using $\tau _{k}$ as the origin for the transformation --- we obtain
a walk which may or may not be self-avoiding.

\smallskip

\textit{Step} 5. If the walk obtained at Step 4 is self-avoiding (in this
case, it belongs to the set $K_{\left( \tau _{0},\tau _{1}\right) }$
(because $\tau =\left( \tau _{0},\tau _{1},...,\tau _{N}\right) $ and $k\geq
1$)), choose it as the new walk, and go to Step 2. If the walk obtained at
Step 4 is not self-avoiding, remove/reject it and keep $\tau ,$ and go to
Step 2.

\smallskip

We call the above algorithm the pivot$^{\text{+}}$ algorithm. This algorithm
is also a Markov chain (a finite one, but it is nonhomogeneous), so, we will
use ``the pivot$^{\text{+}}$ algorithm'' or ``the pivot$^{\text{+}}$ chain''
to denote it. Pivot$^{\text{+}}$ comes from the fact that in each iteration
of the pivot$^{\text{+}}$ algorithm the pivot is positive and --- this also
counts --- has a key role and, on the other hand, from the fact that the
pivot$^{\text{+}}$ algorithm/chain was derived from the pivot
algorithm/chain and is --- we believe this, we have certain reasons for this
--- better than it.

\smallskip

Denote by $p_{n}$ the probability distribution at time $n$ of pivot$^{\text{+%
}}$ chain, $\forall n\geq 0.$ $p_{0}$ is also called the \textit{initial
probability distribution} (of pivot$^{\text{+}}$ chain).

\smallskip

Below we give the fundamental result on the pivot$^{\text{+}}$ chain.

\smallskip

\textbf{THEOREM 2.8.} (i) \textit{The initial probability distribution of
pivot}$^{\text{+}}$\textit{\ chain is} 
\[
p_{0}=\left( \left( p_{0}\right) _{\alpha }\right) _{\alpha \in \Omega _{N}},%
\text{ }\left( p_{0}\right) _{\alpha }=\left\{ 
\begin{array}{cc}
1 & \text{\textit{if} }\alpha =\omega , \\ 
0 & \text{\textit{if} }\alpha \neq \omega ,
\end{array}
\right. 
\]

\noindent $\forall \alpha \in \Omega _{N}$ --- $p_{0}$ \textit{corresponds
to Step} 0 (\textit{of the pivot}$^{\text{+}}$\textit{\ algorithm}).

(ii) \textit{Consider that the self-avoiding walks from }$K_{\left(
0,e_{1}\right) }$\textit{\ are the first }$\left| K_{\left( 0,e_{1}\right)
}\right| $\textit{\ elements }(\textit{self-avoiding walks}) \textit{of }$%
\Omega _{N},$\textit{\ the self-avoiding walks from }$K_{\left(
0,e_{2}\right) }$\textit{\ are the second }$\left| K_{\left( 0,e_{2}\right)
}\right| $\textit{\ elements of }$\Omega _{N},$\textit{\ ..., the
self-avoiding walks from }$K_{\left( 0,e_{-d}\right) }$\textit{\ are the
last }$\left| K_{\left( 0,e_{-d}\right) }\right| $\textit{\ elements of }$%
\Omega _{N}.$\textit{\ The transition matrices of pivot}$^{\text{+}}$\textit{%
\ chain are }(\textit{this chain is nonhomogeneous})\textit{, say, }$P_{1},$%
\textit{\ }$P_{2},$\textit{\ }$P_{2},$\textit{\ }$...$\textit{\ --- }$P_{1}$%
\textit{\ corresponds to Step }1, \textit{and is also the first matrix of
reference chain};\textit{\ }$P_{2},$\textit{\ }$P_{2},$\textit{\ }$...$%
\textit{\ correspond to Steps} 2-5 ---, \textit{where} (\textit{see Step} 1) 
\[
P_{1}=\left( \left( P_{1}\right) _{\alpha \beta }\right) _{\alpha ,\beta \in
\Omega _{N}},\text{ }\left( P_{1}\right) _{\alpha \beta }=\left\{ 
\begin{array}{ll}
\frac{1}{2d}\smallskip & \text{if }\beta \in B_{N}, \\ 
0 & \text{if }\beta \notin B_{N},
\end{array}
\right. 
\]

\noindent $\forall \alpha ,\beta \in \Omega _{N},$ \textit{and }(\textit{see
Steps} 2-5) $P_{2}$\textit{\ is a block diagonal matrix}, 
\[
P_{2}=\left( 
\begin{array}{cccccccc}
Q_{e_{1}} &  &  &  &  &  &  &  \\ 
& Q_{e_{2}} &  &  &  &  &  &  \\ 
&  & \ddots &  &  &  &  &  \\ 
&  &  & Q_{e_{d}} &  &  &  &  \\ 
&  &  &  & Q_{-e_{1}} &  &  &  \\ 
&  &  &  &  & Q_{-e_{2}} &  &  \\ 
&  &  &  &  &  & \ddots &  \\ 
&  &  &  &  &  &  & Q_{-e_{d}}
\end{array}
\right) , 
\]
\[
Q_{\gamma _{1}}=\left( \left( Q_{\gamma _{1}}\right) _{\alpha \beta }\right)
_{\alpha ,\beta \in K_{\left( 0,\gamma _{1}\right) }},\text{ }\forall \alpha
,\beta \in K_{\left( 0,\gamma _{1}\right) },\text{ }\forall \gamma _{1}\in
E_{d} 
\]
\noindent (\textit{for the computation of entries of stochastic matrices }$%
Q_{\gamma _{1}},$\textit{\ }$\gamma _{1}\in E_{d},$\textit{\ when this thing
is necessary and possible, use Steps} 2-5).

(iii) $\rho _{\gamma _{1}}Q_{\gamma _{1}}=\rho _{\gamma _{1}},$\textit{\ }$%
\forall \gamma _{1}\in E_{d},$\textit{\ where }$\rho _{\gamma _{1}}$\textit{%
\ is the uniform probability distribution on }$K_{\left( 0,\gamma
_{1}\right) },$\textit{\ }$\forall \gamma _{1}\in E_{d}.$

(iv) $Q_{\gamma _{1}}^{n}\longrightarrow e^{\prime }\rho _{\gamma _{1}}$ 
\textit{as} $n\longrightarrow \infty ,$\textit{\ }$\forall \gamma _{1}\in
E_{d}$ (\textit{here, }$e$\textit{\ has }$\left| K_{\left( 0,\gamma
_{1}\right) }\right| $\textit{\ components},\textit{\ }$\forall \gamma
_{1}\in E_{d}$).

(v) $P_{1}P_{2}^{n}\longrightarrow e^{\prime }\pi $ \textit{as} $%
n\longrightarrow \infty $ (\textit{recall that }$\pi $\textit{\ is the
uniform probability distribution on }$\Omega _{N};$ \textit{here, }$e$%
\textit{\ has }$\left| \Omega _{N}\right| $\textit{\ components}).

(vi) $p_{n}\longrightarrow \pi $ \textit{as} $n\longrightarrow \infty $ (%
\textit{therefore, }$\pi $\textit{\ is the limit probability distribution of
pivot}$^{\text{+}}$\textit{\ chain}).

\smallskip

\textit{Proof}. (i) Obvious.

(ii) Obvious --- on $P_{1},$ no problem; as to $P_{2},$ when the chain
enters a set $K_{\left( 0,\gamma _{1}\right) }$ for some $\gamma _{1}\in
E_{d}$, see Step 1, it stays forever in this set, see Steps 2-5.

(iii) This follows from 
\[
\left( \rho _{\gamma _{1}}\right) _{\alpha }\left( Q_{\gamma _{1}}\right)
_{\alpha \beta }=\left( \rho _{\gamma _{1}}\right) _{\beta }\left( Q_{\gamma
_{1}}\right) _{\beta \alpha },\text{ }\forall \alpha ,\beta \in K_{\left(
0,\gamma _{1}\right) },\text{ }\forall \gamma _{1}\in E_{d}, 
\]
\noindent because 
\[
\left( \rho _{\gamma _{1}}\right) _{\alpha }=\left( \rho _{\gamma
_{1}}\right) _{\beta }=\frac{1}{a_{N}},\text{ }\forall \alpha ,\beta \in
K_{\left( 0,\gamma _{1}\right) },\text{ }\forall \gamma _{1}\in E_{d}, 
\]
\noindent and $Q_{\gamma _{1}}$ is symmetric (any orthogonal transformation
is invertible, and see/use Steps 2-5), $\forall \gamma _{1}\in E_{d}.$

(iv) Let $\gamma _{1}\in E_{d}.$ Let $\alpha ,\beta \in K_{\left( 0,\gamma
_{1}\right) }.$ Then 
\[
\left( Q_{\gamma _{1}}\right) _{\alpha \beta }>0\text{ if and only if }%
P_{\alpha \beta }>0\text{ (not }\left( Q_{\gamma _{1}}\right) _{\alpha \beta
}=P_{\alpha \beta }\text{),} 
\]
\noindent $P$ is the transition matrix of pivot chain. This statement is
obvious. By Theorem 2.3(ii), $P_{K_{\left( 0,\gamma _{1}\right)
}}^{K_{\left( 0,\gamma _{1}\right) }}$ is irreducible. Moreover, it is
aperiodic irreducible because, moreover, $P_{\alpha \alpha }>0,$ $\forall
\alpha \in K_{\left( 0,\gamma _{1}\right) }.$ So, $Q_{\gamma _{1}}$ is
aperiodic irreducible. Recall that $Q_{\gamma _{1}}$ is symmetric (see the
proof of (iii)). Further, we have --- see (iii) and Theorem 2.5(iii) and its
proof again --- 
\[
Q_{\gamma _{1}}^{n}\longrightarrow e^{\prime }\rho _{\gamma _{1}}\text{%
\textit{\ }as\textit{\ }}n\longrightarrow \infty . 
\]

(v) We have 
\[
\lim_{n\rightarrow \infty }P_{1}P_{2}^{n}=P_{1}\lim_{n\rightarrow \infty
}P_{2}^{n}= 
\]
\[
=P_{1}\left( 
\begin{array}{cccccccc}
e^{\prime }\rho _{e_{1}} &  &  &  &  &  &  &  \\ 
& e^{\prime }\rho _{e_{2}} &  &  &  &  &  &  \\ 
&  & \ddots &  &  &  &  &  \\ 
&  &  & e^{\prime }\rho _{e_{d}} &  &  &  &  \\ 
&  &  &  & e^{\prime }\rho _{-e_{1}} &  &  &  \\ 
&  &  &  &  & e^{\prime }\rho _{-e_{2}} &  &  \\ 
&  &  &  &  &  & \ddots &  \\ 
&  &  &  &  &  &  & e^{\prime }\rho _{-e_{d}}
\end{array}
\right) . 
\]
\noindent By (ii) we have 
\[
P_{1}^{\left\{ \beta \right\} }=\left( 
\begin{array}{c}
\begin{array}{c}
\frac{1}{2d}\smallskip \\ 
\frac{1}{2d}\smallskip \\ 
\vdots \smallskip \\ 
\frac{1}{2d}
\end{array}
\end{array}
\right) 
\]
\noindent ($P_{1}^{\left\{ \beta \right\} }$ is the column $\beta $ of $%
P_{1} $, see Section 1 (see the definition of matrix $P_{U}^{V}...$)) and 
\[
P_{1}^{K_{\left( 0,\gamma _{1}\right) }-\left\{ \beta \right\} }=0\text{ (a
zero matrix), }\forall \gamma _{1}\in E_{d},\text{ }\forall \beta \in
B_{N}\cap K_{\left( 0,\gamma _{1}\right) } 
\]
\noindent ($B_{N}\cap K_{\left( 0,\gamma _{1}\right) }$ has only one element
--- the straight self-avoiding walk from $K_{\left( 0,\gamma _{1}\right) }$%
). By Theorem 2.6, 
\[
\left( K_{\left( 0,\gamma _{1}\right) }\right) _{\left( 0,\gamma _{1}\right)
\in \Omega _{1}} 
\]
\noindent is a partition of $\Omega _{N}.$ Set 
\[
\Delta _{1}=\left( \Omega _{N}\right) ,\text{ }\Delta _{2}=\left( K_{\left(
0,\gamma _{1}\right) }\right) _{\left( 0,\gamma _{1}\right) \in \Omega _{1}},%
\text{ }\Delta _{3}=\left( \left\{ \alpha \right\} \right) _{\alpha \in
\Omega _{N}}. 
\]
\noindent We have 
\[
P_{1}\in G_{\Delta _{1},\Delta _{2}}\text{ and }\lim_{n\rightarrow \infty
}P_{2}^{n}\in G_{\Delta _{2},\Delta _{3}}. 
\]
\noindent Further, by Theorem 1.6, 
\[
P_{1}\lim_{n\rightarrow \infty }P_{2}^{n} 
\]
\noindent is a stable matrix and 
\[
\left( P_{1}\lim_{n\rightarrow \infty }P_{2}^{n}\right) _{\left\{ \alpha
\right\} }=P_{1}^{-+}\left( \lim_{n\rightarrow \infty }P_{2}^{n}\right)
^{-+}= 
\]
\[
=\left( 
\begin{array}{cccc}
\frac{1}{2d} & \frac{1}{2d} & \cdots & \frac{1}{2d}
\end{array}
\right) \left( 
\begin{array}{cccccccc}
\rho _{e_{1}} &  &  &  &  &  &  &  \\ 
& \rho _{e_{2}} &  &  &  &  &  &  \\ 
&  & \ddots &  &  &  &  &  \\ 
&  &  & \rho _{e_{d}} &  &  &  &  \\ 
&  &  &  & \rho _{-e_{1}} &  &  &  \\ 
&  &  &  &  & \rho _{-e_{2}} &  &  \\ 
&  &  &  &  &  & \ddots &  \\ 
&  &  &  &  &  &  & \rho _{-e_{d}}
\end{array}
\right) = 
\]
\[
=\left( 
\begin{array}{cccc}
\frac{1}{2da_{N}} & \frac{1}{2da_{N}} & \cdots & \frac{1}{2da_{N}}
\end{array}
\right) =\left( 
\begin{array}{cccc}
\frac{1}{c_{N}} & \frac{1}{c_{N}} & \cdots & \frac{1}{c_{N}}
\end{array}
\right) =\pi ,\text{ }\forall \alpha \in \Omega _{N} 
\]
\noindent ($\rho _{\gamma _{1}}=\left( 
\begin{array}{cccc}
\frac{1}{a_{N}} & \frac{1}{a_{N}} & \cdots & \frac{1}{a_{N}}
\end{array}
\right) ,$ $\forall \gamma _{1}\in E_{d},$ and, see Theorem 2.7, $%
c_{N}=2da_{N}$). So, 
\[
\lim_{n\rightarrow \infty }P_{1}P_{2}^{n}=e^{\prime }\pi . 
\]

(vi) We have 
\[
p_{n}=p_{0}P_{1}P_{2}^{n-1},\forall n\geq 1. 
\]
\noindent By (v), 
\[
\lim_{n\rightarrow \infty }p_{n}=\lim_{n\rightarrow \infty
}p_{0}P_{1}P_{2}^{n-1}=p_{0}\lim_{n\rightarrow \infty
}P_{1}P_{2}^{n-1}=p_{0}e^{\prime }\pi =\pi .\text{ }\square 
\]

\smallskip

We comment on Theorem 2.8 and its proof.

\smallskip

\textbf{1.} $P_{1}$ and $\lim\limits_{n\rightarrow \infty }P_{2}^{n}$ (for
these matrices, see (ii) and the proof of (v)) are the matrices of reference
chain associated with the reference method with Steps I and II. $P_{1},$ $%
e^{\prime }\pi \in G_{\Delta _{1},\Delta _{2}}$ (see also Remark 2.1(f) in
[13]) and (an interesting fact) $P_{1}\backsim e^{\prime }\pi ;$ $%
\lim\limits_{n\rightarrow \infty }P_{2}^{n},$ $e^{\prime }\pi \in G_{\Delta
_{2},\Delta _{3}},$ but $\lim\limits_{n\rightarrow \infty }P_{2}^{n}\not%
{\backsim}e^{\prime }\pi .$ The construction of $P_{1}$ is due to the
reference chain. $P_{1}$ can be replaced with any matrix which is similar to
it (equivalently, to $e^{\prime }\pi $) with respect to $\Delta _{1}$ and $%
\Delta _{2}$, see Theorem 1.9, in this case, Step 1 (of the pivot$^{\text{+}%
} $ algorithm) must be modified accordingly --- but $P_{1}$ is very good; it
was optimally constructed, \textit{i.e.}, it was constructed to have the
smallest possible number of positive entries (see also Examples 1.8 and
1.10).

\smallskip

\textbf{2.} The product $P_{1}\lim\limits_{n\rightarrow \infty }P_{2}^{n}$
(see (v) and its proof) can be computed directly, without Theorem 1.6,
because it is easy to compute it , but we preferred to use this theorem ---
and thus we gave a new application of this theorem!

\smallskip

\textbf{3.} When the chain enters a set $K_{\left( 0,\gamma _{1}\right) }$
for some $\gamma _{1}\in E_{d}$ --- this thing happens at time $1$ (not at
time $2$) ---, it stays forever in this set (see (ii) and its proof). So, at
times $1,2,...$ we work with a chain having $K_{\left( 0,\gamma _{1}\right)
} $ as the state space, a state space $2d$ times smaller than $\Omega _{N},$
see Theorem 2.7 --- \textit{e.g}., for $d=2$ and $N=10,$ we have $\left|
\Omega _{N}\right| =44100,$ see, \textit{e.g.}, [9, p. 394], and $\left|
K_{\left( 0,\gamma _{1}\right) }\right| =11025.$

\smallskip

\textbf{4.} $\pi $ is a wavy probability distribution with respect to the
partitions 
\[
\Delta _{1}=\left( \Omega _{N}\right) ,\text{ }\Delta _{2}=\left( K_{\left(
0,\gamma _{1}\right) }\right) _{\left( 0,\gamma _{1}\right) \in \Omega _{1}},%
\text{ and }\Delta _{3}=\left( \left\{ \alpha \right\} \right) _{\alpha \in
\Omega _{N}} 
\]

\noindent --- any order relation on $\Omega _{N}$ can be used because $\pi $
is uniform; we have $c_{N}!$ order relations on $\Omega _{N}$. (For the wavy
probability distributions, see, \textit{e.g.}, [16]--[17]; the wavy
probability distributions with respect to three or more partitions can lead
to good or very good results, see [17, p. 215].)

\smallskip

We also cannot run the pivot$^{\text{+}}$ algorithm/chain forever... --- for
the completion, see the first paragraph after Theorem 2.5 and its proof.

\smallskip 

It is worthy to note that the pivot$^{\text{+}}$ algorithm works both when
the observables being measured are invariant with respect to the symmetries
of lattice and when these are not. The pivot$^{\text{+}}$ algorithm without
Step 1 (skipping Step 1) is an algorithm which, obviously, is not better
than the pivot$^{\text{+}}$ algorithm, but can be used when the observables
being measured are invariant with respect to the symmetries of lattice (%
\textit{e.g.}, the end-to-end distance (of a self-avoiding walk) is
invariant with respect to these transformations). This algorithm is known,
see, \textit{e.g.}, the variant 1 of pivot algorithm in [10, p. 117] (1988)
and, somehow equivalently, but better (correct statements...), the first
paragraph before Theorem 9.4.4 in [9, p. 324] (1996 --- eight years later).

\smallskip 

We must compare the pivot$^{\text{+}}$ algorithm with the pivot algorithm
--- 
\[
which\text{ }of\text{ }them\text{ }is\text{ }better? 
\]

\smallskip

\noindent The pivot$^{\text{+}}$ algorithm samples the pivot from the set $%
\left\langle N-1\right\rangle $ according to the uniform probability on this
set while the pivot algorithm samples the pivot from the set $\left\langle
\left\langle N-1\right\rangle \right\rangle $ according to the uniform
probability on this set, so, on the pivot, the pivot$^{\text{+}}$ algorithm
is better than the pivot algorithm. As to the speed of convergence, we give
a conjecture --- Comment 3 and the previous fact (see the previous sentence)
seem that support this conjecture.

\smallskip

Set 
\[
\left| \left| x\right| \right| _{1}=\sum\limits_{i=1}^{n}\left| x_{i}\right|
, 
\]

\noindent where $x=\left( x_{1},x_{2},...,x_{n}\right) \in \Bbb{R}^{n}$ and $%
n\geq 1.$ $\left| \left| x\right| \right| _{1}$ is called the 1-\textit{norm}
\textit{of }$x$; $\left| \left| \cdot \right| \right| _{1}$ is called the 1-%
\textit{norm.}

\smallskip

\textit{Conjecture} 2.9. The pivot$^{\text{+}}$ algorithm/chain is faster
than the pivot algorithm/chain, \textit{i.e.}, 
\[
\left| \left| p_{n}-\pi \right| \right| _{1}\leq \left| \left| q_{n}-\pi
\right| \right| _{1},\text{ }\forall n\geq n_{0} 
\]
\noindent where $n_{0}\geq 1$ (does $n_{0}$ depend on $N$? can we replace ``$%
\forall n\geq n_{0},$ where $n_{0}\geq 1$'' with ``$\forall n\geq 2$''?).

\smallskip

Concerning the above conjecture, obviously, the technology is not taken into
account.

\smallskip

Generally speaking, the comparison of (finite or not) Markov chains, when
this makes sense, on the speed of convergence is or seems to be a very
hard/difficult problem, and is important --- \textit{e.g.}, some of the
known Markov chains for sampling are not fast or are not too fast, so, we
must find, if possible, faster ones.

\smallskip

We consider that the homogeneous finite Markov chain framework is too narrow
to design very good Markov chains for sampling and, besides sampling, if
possible, for other things (Markov chain Monte Carlo methods or, best, see, 
\textit{e.g.}, [16], Markov chain methods) --- the nonhomogeneous finite
Markov chain framework is better.

\smallskip 

Our aim was twofold: to give new applications of the $G$ method and to
construct algorithms/chains for the self-avoiding walk better than the pivot
algorithm.
\newpage
\bigskip

\begin{center}
\textbf{REFERENCES}
\end{center}

\bigskip \ 

[1] N. Clisby, \textit{Enumerative combinatorics of lattice polymers}.
Notices Amer. Math. Soc. \textbf{68} (2021), 504$-$515.

[2] P.J. Flory, \textit{Principles of Polymer Chemistry}. Cornell University
Press, Ithaca, 1953.

[3] J.M. Hammersley and D.C. Handscomb, \textit{Monte Carlo Methods}.
Methuen, London \& Wiley, New York, 1964.

[4] M. Iosifescu, \textit{Finite Markov Processes and Their Applications}.
Wiley, Chichester \& Ed. Tehnic\u {a}, Bucharest, 1980; corrected
republication by Dover, Mineola, N.Y., 2007.

[5] D.L. Isaacson and R.W. Madsen, \textit{Markov Chains}:\textit{\ Theory
and Applications.} Wiley, New York, 1976; republication by Krieger, 1985.

[6] M. Lal, `\textit{Monte Carlo}' \textit{computer simulations of chain
molecules}. I. Mol. Phys. \textbf{17} (1969), 57$-$64.

[7] G.F. Lawler and L.N. Coyle, \textit{Lectures on Contemporary Probability}%
. AMS, Providence, Rhode Island, 1999.

[8] N. Madras, \textit{Lectures on Monte Carlo Methods}. AMS, Providence,
Rhode Island, 2002.

[9] N. Madras and G. Slade, \textit{The Self-Avoiding Walk}. Birkh\"{a}user,
Boston, 1996.

[10] N. Madras and A.D. Sokal, \textit{The pivot algorithm}:\textit{\ a
highly efficient Monte Carlo method for the self-avoiding walk.} J. Stat.
Phys. \textbf{50} (1988), 109--186.

[11] U. P\u {a}un, \textit{New classes of ergodicity coefficients, and
applications}. Math. Rep. (Bucur.) \textbf{6}(\textbf{56}) (2004), 141$-$158.

[12] U. P\u {a}un, \textit{Ergodicity coefficients of several matrices}:%
\textit{\ new results and applications}. Rev. Roumaine Math. Pures Appl.%
\textbf{\ 55 }(2010), 53$-$77.

[13] U. P\u {a}un, $G_{\Delta _{1},\Delta _{2}}$ \textit{in action}. Rev.
Roumaine Math. Pures Appl.\textbf{\ 55 }(2010), 387$-$406.

[14] U. P\u {a}un, \textit{A hybrid Metropolis-Hastings chain}. Rev.
Roumaine Math. Pures Appl.\textbf{\ 56 }(2011), 207$-$228.

[15] U. P\u {a}un, $G$\textit{\ method in action}: \textit{from exact
sampling to approximate one}. Rev. Roumaine Math. Pures Appl.\textbf{\ 62 }%
(2017), 413$-$452.

[16] U. P\u {a}un, \textit{Ewens distribution on }$\Bbb{S}_{n}$ \textit{is a
wavy probability distribution with respect to }$n$\textit{\ partitions.} An.
Univ. Craiova Ser. Mat. Inform. \textbf{47} (2020), 1$-$24.

[17] U. P\u {a}un, $\Delta $-\textit{wavy probability distributions and
Potts model.} An. Univ. Craiova Ser. Mat. Inform. \textbf{49} (2022), 208$-$%
249.

[18] U. P\u {a}un, $G^{+}$ \textit{method in action}:\textit{\ new classes
of nonnegative matrices}, \textit{with results.} arXiv:2304.05227.

[19] E. Seneta, \textit{Non-negative Matrices and Markov Chains}, 2nd
Edition. Springer-Verlag, Berlin, 1981; revised printing, 2006.

[20] C. Vanderzande, \textit{Lattice Models of Polymers}. Cambridge
University Press, Cambridge, 1998.

\bigskip \ 
\[
\begin{array}{ccc}
\mathit{August}\text{ }\mathit{27,}\text{ }\mathit{2024} &  & \mathit{%
Romanian\ Academy} \\ 
&  & \mathit{Gheorghe\ Mihoc}-\mathit{Caius\ Iacob}\text{ }\mathit{Institute}
\\ 
&  & \mathit{of\ Mathematical\ Statistics}\text{ }\mathit{and\ Applied\
Mathematics} \\ 
&  & \mathit{Calea\ 13\ Septembrie\ nr.\ 13} \\ 
&  & \mathit{050711\ Bucharest\ 5,\ Romania} \\ 
&  & \mathit{upterra@gmail.com}
\end{array}
\]

\end{document}